## UM NOVO MÉTODO ANALÍTICO PARA RESOLVER A EQUAÇÃO DE TELEGRAFIA LINEAR E NÃO-LINEAR

## A NOVEL ANALYTIC METHOD FOR SOLVING LINEAR AND NON-LINEAR TELEGRAPH EQUATION

Al-JABERI, Ahmed K.[1], HAMEED, Ehsan M.[2], Abdul-WAHAB, Mohammed S.[3]

[1,3] Department of Mathematics, College of Education for Pure Sciences, University of Basrah, Basrah, Iraq.

[2] Department of Mathematics, College of Computer Science and Mathematics, University of Thi-Qar, Thi-Qar, Iraq.

*\* Corresponding author*
*e-mail: ahmed.shanan@uobasrah.edu.iq*

### RESUMO

A modelagem de muitos fenômenos em vários campos, como matemática, física, química, engenharia, biologia e astronomia, é feita pelas equações diferenciais parciais não lineares (PDE). A equação do telégrafo hiperbólico é uma delas, onde descreve as vibrações de estruturas (por exemplo, edifícios, vigas e máquinas) e é a base para equações fundamentais da física atômica. Existem vários métodos analíticos e numéricos para resolver a equação do telégrafo. Uma solução analítica considera enquadrar o problema de uma forma bem compreendida e calcular a resolução exata. Também ajuda a entender as respostas para o problema em termos de precisão e convergência. Esses métodos analíticos têm limitações com precisão e convergência. Portanto, um novo método analítico aproximado é proposto para lidar com restrições neste artigo. Este método usa as séries de Taylors em sua derivação. O método proposto foi usado para resolver a equação hiperbólica de segunda ordem (equação Telegraph) com a condição inicial. Três exemplos foram apresentados para verificar a eficácia, precisão e convergência do método. As soluções do método proposto também foram comparadas com as obtidas pelo método de decomposição adomiana (ADM) e pelo método de análise de homotopia (HAM). A técnica é fácil de implementar e produz resultados precisos. Em particular, esses resultados mostram que o método proposto é eficiente e melhor que os outros métodos em termos de precisão e convergência

**Palavras-chave**: *Equação do Telégrafo Hiperbólico, Série de Taylor, Solução Analítica, Precisão, Operador Não Linear.*

### ABSTRACT

The modeling of many phenomena in various fields such as mathematics, physics, chemistry, engineering, biology, and astronomy is done by the nonlinear partial differential equations (PDE). The hyperbolic telegraph equation is one of them, where it describes the vibrations of structures (e.g., buildings, beams, and machines) and are the basis for fundamental equations of atomic physics. There are several analytical and numerical methods are used to solve the telegraph equation. An analytical solution considers framing the problem in a well-understood form and calculating the exact resolution. It also helps to understand the answers to the problem in terms of accuracy and convergence. These analytic methods have limitations with accuracy and convergence. Therefore, a novel analytic approximate method is proposed to deal with constraints in this paper. This method uses the Taylors' series in its derivation. The proposed method has used for solving the second-order, hyperbolic equation (Telegraph equation) with the initial condition. Three examples have been presented to check the effectiveness, accuracy, and convergence of the method. The solutions of the proposed method also compared with those obtained by the Adomian decomposition method (ADM), and the Homotopy analysis method (HAM). The technique is easy to implement and produces accurate results. In particular, these results display that the proposed method is efficient and better than the other methods in terms of accuracy and convergence.

## 1. INTRODUCTION:

Telegraph's equation is one of the most investigated problems in over four centuries. Generally, it uses to model the vibrations of structures (e.g., buildings, beams, and machines) and are the basis for fundamental equations of atomic physics (Kolesov and Rozov 2000) and (Alonso, Mawhin, and Ortega 1999).

Linear hyperbolic telegraph equation with constant coefficients forms a blend between wave propagation and diffusion by inserting a term that accounts for impacts of finite velocity to the standard heat or mass transport equation (El-Azab and El-Gamel 2007). In general, the hyperbolic telegraph equation is used in signal analysis for transmission and propagation of electrical signals (A. C. Metaxas 1983) and also has applications in other fields (see Roussy and Pearce 1998) and the references therein).

Recently, several numerical methods have been developed to solve hyperbolic partial differential equations (PDEs), such as alternating direction implicit Method (Mohanty, Jain, and Arora 2002), Chebyshev tau method (Saadatmandi and Dehghan 2010), Interpolating Scaling Functions Method (Lakestani and Saray 2010), Radial Basis Functions method (De Su, Jiang, and Jiang 2013), collocation points and approximating the solution using thin-plate splines radial basis function method (Dehghan and Shokri 2008), the modified cubic B-Spline differential quadrature method (Singh 2016), a collocation method (Bhatia 2014), reproducing kernel method (Mustafa Inc 2013), finite difference method (Mohanty, Jain, and George 1996), and (Gao and Chi 2007), finite volume method (Sheng and Zhang 2018), Galerkin method (Liu and Gu 2004), improved element-free Galerkin method (Zhang, Deng, and Liew 2014), and differential quadrature method (Mittal and Jiwari 2009).

These techniques consider a necessitate as computational resources to solve some problems that appear in other sciences such as in image processing. Still, these techniques perhaps complicated and require a high computational cost, which means they consume time and effort to obtain accurate solutions for different nonlinear PDEs. Furthermore, applying these techniques may be needed for transforming. Thus, they do not avoid linearization, discretization, or any realistic assumption to providing an efficient solution. Therefore, it is considered positive points for analytical methods to find the most efficient and high accurate approximate and exact solutions for the linear and nonlinear differential equations.

There are many used methods to find the exact solution, such as ($G'/G$ ) expansion method (Naher, Hasibun 2012), sine-cosine method (Yusufoğlu and Bekir 2006), homogeneous balance method (Eslami, Fathi Vajargah, and Mirzazadeh 2014), etc. Also, there are other types of analytical techniques that combine between exact and approximate methods called semi-analytical methods. For instance, Adomian decomposition method (ADM), (Adomian 1994) and (Biazar and Ebrahimi 2007), the extended of ADM (Al-Badrani et al. 2016), Differential transform method (DTM) (Zhou 1986), Variational Iteration Method (VIM) (Biazar, Ebrahimi, and Ayati 2008), Reduced differential transform method (Mukesh K. Awasthi, R. K. Chaurasia 2013) and (Srivastava, Awasthi, and Chaurasia 2014), Daftardar-Gejji-Jafaris (DGJ) method (Sari, Gunay, and Gurarslan 2014), Reimann method (ZAVALANI 2009), Sumudu method (Ayad Ghazi Naser Al-Shammari, Wurood R. Abd AL-Hussein 2018), the modified simple equation method (Taghizadeh et al. 2012), Chebyshev Cardinal functions are used for solution of the second-order one-dimensional telegraph Equation (Dehghan and Lakestani 2009), Wavefront solutions of a nonlinear telegraph equation (Gilding and Kersner 2013), Homotopy analysis method (HAM) (Raftari and Yildirim 2012), Homotopy analysis Sumudu transform method (Rathore et al. 2012) and (Jaradat et al. 2018), and Laplace transform method (LDM) (Jaradat et al. 2018).

One of these methods which has received much concern is the Adomian decomposition method (ADM) (Biazar and Ebrahimi 2007). The ADM has been employed to solve various linear and nonlinear models. The ADM yields a rapidly convergent series solution with much less computational work (Abbaoui and Cherruault 1995). The ADM is unlike the traditional numerical methods, where ADM is used extensively to solve nonlinear differential equations because it works based on calculation Adomian polynomials for non-linear terms (El-Sayed and Kaya 2004) and (Inc 2007).

While, the analytical solutions of differential transform method (DTM) (Zhou 1986), are as a polynomial form which is different from the traditional higher-order Taylor series method because the Taylor series method needs huge computational for large orders. So, the DTM uses a different procedure to obtain an analytic Taylor series solution of the PDE (Abazari and Borhanifar 2010) and (Saeideh Hesam 2011).

On the other hand, many complicated problems in different applied sciences have been successfully solved by the homotopy analysis method (HAM) (Liao 1992), (Kheiri, H. 2011), (Usman et al. 2013) and (Arafa, Rida, and Mohamed 2011). Homotopy perturbation method (HPM) (Raftari and Yildirim 2012), as well has been handled successfully to solve many linear and non-linear PDEs (Roul 2010).

Whereas, the variational iteration method (VIM) (He 1997), can be applied to solve all types of linear or nonlinear differential equations, constant or variable coefficients with homogeneous or inhomogeneous (Yulita Molliq et al. 2009) because it effectively used to solve these nonlinear equations with a good convergent to the exact solutions.

All of these analytical methods are active and strong in finding solutions of linear and non-linear PDEs. Generally, the basic idea of these methods through writing the PDEs in the form of operators after that the inverse operator of time (or space) is taken to calculate the solutions on the interval of the solution domain. In the last years, these analytical methods have been enhanced and modified and to overcome the difficulties encountered in the numerical methods such as finding the exact solution, round off errors, low accuracy, and weak convergence by many researchers. Furthermore, these analytical methods introduced efficiency and high accuracy in finding the analytical exact and approximate solutions for the linear and nonlinear PDEs, which made a good impression of these analytical methods. Therefore, we proposed a new analytical approximate method for solving linear and non-linear telegraph equations.

This study aimed to find approximate analytical solutions to the linear and nonlinear telegraph equation with initial conditions and boundary conditions by using a novel analytical technique, which is considered as extending and developing to that in (wiwatwanich, 2016). This technique is based on the Taylor series, which is efficient to solve nonlinear equations. Also, the survey reveals that no attempt has ever been made to study the current model by using this technique. These reasons stimulated and employed us to solve the linear and intricate nonlinear problems such as the one-dimensional telegraph equation. Several test problems are given, and their results are compared with the solutions obtained by ADM (Biazar and Ebrahimi 2007), and HAM (Raftari and Yildirim 2012) to confirm the excellent accuracy and small absolute errors of the proposed technique.

## 2. MATERIALS AND METHODS:

In this work, it was studied the second-order one-dimensional telegraph equation:

$$\frac{\partial^2 w}{\partial t^2} + 2 \propto \frac{\partial w}{\partial t} + \beta^2 w = \frac{\partial^2 w}{\partial x^2} + h(x,t) + g(w), \quad (1)$$

where $w = w(x,t)$, $\propto$, $\beta$ are real constants, $x$ is distance and $t$ is time, with initial and boundary conditions

$$w(x,0) = h_0(x), \qquad \frac{\partial w}{\partial t}(x,0) = h_1(x), \quad (2)$$

$$w(0,t) = g_0(t), \qquad w(1,t) = g_1(t), \qquad t \geq 0 \quad (3)$$

If $g(w) = 0$, then equation (1) called a linear hyperbolic telegraph equation in a double conductor, the equation (2.1) is satisfied in both the electric voltage and current. For $\alpha > 0$, $\beta = 0$ it represents a damped wave equation, and $\propto > \beta > 0$, which is called a telegraph equation.

### 2.1. Generating an analytical approach

In this section, the basic ideas for constructing a novel analytical approach will be discussed. Let us consider the initial value problems:

$$w_{tt}(x,t) = G(w_t, w, w_x, w_{xx}, \dots), \quad (4)$$

with initial condition

$$w(x,0) = h_0(x), \qquad w_t(x,0) = h_1(x). \quad (5)$$

By using the integral for the two sides of equation (4) from 0 to $t$, we get

$$w_t(x,t) - w_t(x,0) = \int_0^t F[w]\, dt,$$

$$w_t(x,t) - h_1(x) = \int_0^t F[w]\, dt.$$

Then,

$$w_t(x,t) = h_1(x) + \int_0^t F[w]\, dt, \quad (6)$$

where $G[w] = G(w_t, w, w_x, w_{xx}, \dots)$.

Then, when the integral of two sides of equation (6) is used from 0 to $t$, we obtain

$$w(x,t) - w(x,0) = h_1(x)t + \iint_0^t G[w]\, dt\, dt,$$

$$w(x,t) - h_0(x) = h_1(x)t + \iint_0^t G[w]\, dt\, dt.$$

Thus,

$$w(x,t) = h_0(x) + h_1(x)t + \iint_0^t G[w]\, dt\, dt. \tag{7}$$

The Taylor series is extended for $G[w]$ about $t = 0$, which is

$$G[w] = G[w_0] + G'[w_0]t + G''[w_0]\frac{t^2}{2!} + G'''[w_0]\frac{t^3}{3!} + \cdots + G^{(n)}[w_0]\frac{t^n}{n!} + \cdots \tag{8}$$

Substituting equation (8) by equation (7), we get

$$w(x,t) = h_0(x) + h_1(x)t + G[w_0]\frac{t^2}{2!} + G'[w_0]\frac{t^3}{3!} + G''[w_0]\frac{t^4}{4!} + \cdots + G^{(n-2)}[w_0]\frac{t^n}{n!} + \cdots,$$

$$= a_0 + a_1 t + a_2 \frac{t^2}{2!} + a_3 \frac{t^3}{3!} + \cdots + a_n \frac{t^n}{n!} + \cdots, \tag{9}$$

where

$$a_0 = h_0(x),$$

$$a_1 = h_1(x),$$
$$a_2 = G[w_0],$$

$$a_3 = G'[w_0],$$
$$\cdot$$
$$\cdot$$
$$a_n = G^{(n-2)}[w_0].$$

where $n$ is the highest derivative of $u$. The formal of Equation (9) is expand Taylor's series for $w$ about $t = 0$. This means

$$a_0 = w(x,0),$$
$$a_1 = \frac{\partial}{\partial t} w(x,0),$$
$$a_2 = \frac{\partial^2}{\partial t^2} w(x,0),$$

$$a_3 = \frac{\partial^3}{\partial t^3} w(x,0),$$
$$\cdot$$

$$a_n = \frac{\partial^n}{\partial t^n} w(x,0).$$

### 2.1.1 Test Problems

**Example 1.** Solve the following linear telegraph equation (Lakestani and Saray 2010):

$$\frac{\partial^2 w(x,t)}{\partial t^2} + 2 \propto \frac{\partial w(x,t)}{\partial t} + \beta^2 w(x,t) = \frac{\partial^2 w(x,t)}{\partial x^2} + (3 - 4 \propto + \beta^2)e^{-2t}\sinh(x),$$

$$\frac{\partial^2 w(x,t)}{\partial t^2} = -2 \propto \frac{\partial w(x,t)}{\partial t} - \beta^2 w(x,t) + \frac{\partial^2 w(x,t)}{\partial x^2} + (3 - 4 \propto + \beta^2)e^{-2t}\sinh(x). \tag{10}$$

**Solution:**

By the following equation (4), we can note after rewrite equation (10):

$$F[w] = -2 \propto \frac{\partial w(x,t)}{\partial t} - \beta^2 w(x,t) + \frac{\partial^2 w(x,t)}{\partial x^2} + (3 - 4 \propto + \beta^2)e^{-2t}\sinh(x),$$

$$a_0 = h_0(x) = \sinh(x),$$

$$a_1 = h_1(x) = -2\sinh(x),$$

$$a_2 = G[w_0] = -2 \propto \frac{\partial w(x,0)}{\partial t} - \beta^2 w(x,0) + \frac{\partial^2 w(x,0)}{\partial x^2} + (3 - 4 \propto + \beta^2)e^{-2t}\sinh(x)$$

$$= -2 \propto a_1 - \beta^2 a_0 + (a_0)_{xx} + (3 - 4 \propto + \beta^2)e^{-2(0)}\sinh(x)$$

$$= 4 \propto \sinh(x) - \beta^2\sinh(x) + \sinh(x) + (3 - 4 \propto + \beta^2)\sinh(x)$$

$$= (4 \propto -\beta^2)\sinh(x) + \sinh(x) + 3\sinh(x) + (-4 \propto + \beta^2)\sinh(x)$$

$$= 4\sinh(x)$$

$$G'[w] = -2 \propto \frac{\partial^2 w(x,t)}{\partial t^2} - \beta^2\frac{\partial w(x,t)}{\partial t} + \frac{\partial^3 w(x,t)}{\partial t \partial x^2} - 2(3 - 4 \propto + \beta^2)e^{-2t}\sinh(x)$$

$$a_3 = G'[w_0] = -2 \propto \frac{\partial^2 w(x,0)}{\partial t^2} - \beta^2\frac{\partial w(x,0)}{\partial t} + \frac{\partial^3 w(x,0)}{\partial t \partial x^2} - 2(3 - 4 \propto + \beta^2)e^{-2(0)}\sinh(x)$$

$$= -2 \propto a_2 - \beta^2 a_1 + (a_1)_{xx} - 2(3 - 4 \propto + \beta^2)\sinh(x)$$

$$= -8 \propto \sinh(x) + 2\beta^2\sinh(x) - 2\sinh(x) - 2(3 - 4 \propto + \beta^2)\sinh(x)$$

$$= (-8 \propto + 2\beta^2) \sinh(x) - 2\sinh(x) - 6\sinh(x) + (8 \propto - 2\beta^2) \sinh(x)$$

$$= -8 \sinh(x).$$

Now by using equation (9),

$$w(x,t) = a_0 + a_1 t + a_2 \frac{t^2}{2!} + a_3 \frac{t^3}{3!} + \cdots + a_n \frac{t^n}{n!} + \cdots$$

$$= \sinh(x) - 2\sinh(x)\, t + 4\sinh(x)\frac{t^2}{2!} - 8\sinh(x)\frac{t^3}{3!} + \cdots$$

$$= \sinh(x)\left(1 - (2t) + \frac{(2t)^2}{2!} - \frac{(2t)^3}{3!} + \cdots\right)$$

$$= \sinh(x)\, e^{-2t}.$$

Therefore, the graph of exact, analytical solutions, and absolute errors of example 1 for $t = 1$ are given in Figure 1. **Erro! Fonte de referência não encontrada.** shows the absolute error using the proposed technique, ADM, and HAM of example 1 with $t = 1, k = 0.001$ and different values of $\alpha = 10$, and $\beta = 5$.

When the value of $n$ increases, the approximate solution gradually approaches the exact solution and the absolute error will decrease. Also, the values of $\propto$ and $\beta$ do not affect the solution.

**Example 2.** Solve the following linear telegraph equation (Lakestani and Saray 2010):

$$\frac{\partial^2 w(x,t)}{\partial t^2} + 2 \propto \frac{\partial w(x,t)}{\partial t} + \beta^2 w(x,t) = \frac{\partial^2 w(x,t)}{\partial x^2} - 2 \propto \sin(t)\sin(x) + \beta^2 \cos(t)\sin(x),$$

$$\frac{\partial^2 w(x,t)}{\partial t^2} = -2 \propto \frac{\partial w(x,t)}{\partial t} - \beta^2 w(x,t) + \frac{\partial^2 w(x,t)}{\partial x^2} - 2 \propto \sin(t)\sin(x) + \beta^2 \cos(t)\sin(x) \quad (11)$$

**Solution:** By the following equation (4), we can note after rewrite equation (11):

$$G[w] = -2 \propto \frac{\partial w(x,t)}{\partial t} - \beta^2 w(x,t) + \frac{\partial^2 w(x,t)}{\partial x^2} - 2 \propto \sin(t)\sin(x) + \beta^2 \cos(t)\sin(x),$$

$$a_0 = h_0(x) = \sin(x),$$

$$a_1 = h_1(x) = 0,$$

$$a_2 = G[w_0] = -2 \propto \frac{\partial w(x,0)}{\partial t} - \beta^2 w(x,0) + \frac{\partial^2 w(x,0)}{\partial x^2} - 2 \propto \sin(0)\sin(x) + \beta^2 \cos(0)\sin(x)$$

$$= -2 \propto a_1 - \beta^2 a_0 + (a_0)_{xx} + \beta^2 \sin(x)$$

$$= -2 \propto (0) - \beta^2 \sin(x) + (\sin(x))_{xx} + \beta^2 \sin(x)$$

$$= -\sin(x),$$

$$G'[w] = -2 \propto \frac{\partial^2 w(x,t)}{\partial t^2} - \beta^2 \frac{\partial w(x,t)}{\partial t} + \frac{\partial^3 w(x,t)}{\partial t \partial x^2} - 2 \propto \cos(t)\sin(x) - \beta^2 \sin(t)\sin(x),$$

$$a_3 = G'[w_0] = -2 \propto \frac{\partial^2 w(x,0)}{\partial t^2} - \beta^2 \frac{\partial w(x,0)}{\partial t} + \frac{\partial^3 w(x,0)}{\partial t \partial x^2} - 2 \propto \cos(0)\sin(x) - \beta^2 \sin(0)\sin(x)$$

$$= -2 \propto a_2 - \beta^2 a_1 + (a_1)_{xx} - 2 \propto \sin(x)$$

$$= 2 \propto \sin(x) - \beta^2(0) + (0)_{xx} - 2 \propto \sin(x)$$

$$= 0,$$

$$G''[w] = -2 \propto \frac{\partial^3 w(x,t)}{\partial t^3} - \beta^2 \frac{\partial^2 w(x,t)}{\partial t^2} + \frac{\partial^4 w(x,t)}{\partial t^2 \partial x^2} + 2 \propto \sin(t)\sin(x) - \beta^2 \cos(t)\sin(x),$$

$$a_4 = G'[w_0] = -2 \propto \frac{\partial^3 w(x,0)}{\partial t^3} - \beta^2 \frac{\partial^2 w(x,0)}{\partial t^2} + \frac{\partial^4 w(x,0)}{\partial t^2 \partial x^2} + 2 \propto \sin(0)\sin(x) - \beta^2 \cos(0)\sin(x)$$

$$= -2 \propto a_3 - \beta^2 a_2 + (a_2)_{xx} - \beta^2 \sin(x)$$

$$= -2 \propto (0) + \beta^2 \sin(x) + \sin(x) - \beta^2 \sin(x)$$

$$= \sin(x),$$

Now, by using the equation (9),

$$w(x,t) = a_0 + a_1 t + a_2 \frac{t^2}{2!} + a_3 \frac{t^3}{3!} + \cdots + a_n \frac{t^n}{n!} + \cdots$$

$$= \sin(x) + (0)t + (-\sin(x))\frac{t^2}{2!} + (0)\frac{t^3}{3!} + (\sin(x))\frac{t^4}{4!} + \cdots$$

$$= \sin(x) - \sin(x)\frac{t^2}{2!} + \sin(x)\frac{t^4}{4!} + \cdots$$

$$= \sin(x)\left(1 - \frac{t^2}{2!} + \frac{t^4}{4!} + \cdots\right)$$

$$= \sin(x)\cos(t).$$

Therefore, the graph of exact, analytical solutions, and absolute errors of example 2 for $t = 1$ are given in **Erro! Fonte de referência não encontrada.**. **Erro! Fonte de referência não encontrada.** shows the absolute error using the proposed technique, ADM, and HAM of example 2 with $t = 1, k = 0.001$ and different values of $\alpha = 10$, and $\beta = 5$.

When the value of $n$ increases, the approximate solution gradually approaches the exact solution and the absolute error will decrease. Also, the values of $\propto$ and $\beta$ do not affect the solutions.

**Example 3.** Solve the following nonlinear telegraph equation (Bhatia 2014):

$$\frac{\partial^2 w(x,t)}{\partial t^2} + \frac{\partial w(x,t)}{\partial t} = \frac{\partial^2 w(x,t)}{\partial x^2} - 2\sin(w(x,t)) - \pi^2 e^{-t}\cos(\pi x) + 2\sin(e^{-t}(1-\cos(\pi x))),$$

$$\frac{\partial^2 w(x,t)}{\partial t^2} = -\frac{\partial w(x,t)}{\partial t} + \frac{\partial^2 w(x,t)}{\partial x^2} - 2\sin(w(x,t)) - \pi^2 e^{-t}\cos(\pi x) + 2\sin(e^{-t}(1-\cos(\pi x))).$$

**Solution:** By the following equation (4), we can note after rewrite equation (12):

$$G[w] = -\frac{\partial w(x,t)}{\partial t} + \frac{\partial^2 w(x,t)}{\partial x^2} - 2\sin(w(x,t)) - \pi^2 e^{-t}\cos(\pi x) + 2\sin(e^{-t}(1-\cos(\pi x))),$$

$$a_0 = h_0(x) = 1 - \cos(\pi x),$$

$$a_1 = h_1(x) = -1 + \cos(\pi x),$$

$$a_2 = G[w_0] = -\frac{\partial w(x,0)}{\partial t} + \frac{\partial^2 w(x,0)}{\partial x^2} - 2\sin(w(x,0)) - \pi^2 e^{-(0)}\cos(\pi x) + 2\sin(e^{-(0)}(1-\cos(\pi x)))$$

$$= -a_1 + (a_0)_{xx} - 2\sin(a_0) - \pi^2\cos(\pi x) + 2\sin((1-\cos(\pi x)))$$

$$= 1 - \cos(\pi x) + \pi^2\cos(\pi x) - 2\sin(1-\cos(\pi x)) - \pi^2\cos(\pi x) + 2\sin((1-\cos(\pi x)))$$

$$= 1 - \cos(\pi x),$$

$$G'[w] = -\frac{\partial^2 w(x,t)}{\partial t^2} + \frac{\partial^3 w(x,t)}{\partial t \partial x^2} - 2\cos(w(x,t))\frac{\partial w(x,t)}{\partial t} + \pi^2 e^{-t}\cos(\pi x) + 2\cos(e^{-t}(1-\cos(\pi x)))e^{-t}(-1+\cos(\pi x)),$$

$$a_3 = G'[w_0] = -\frac{\partial^2 w(x,0)}{\partial t^2} + \frac{\partial^3 w(x,0)}{\partial t \partial x^2} - 2\cos(w(x,0))\frac{\partial w(x,0)}{\partial t} + \pi^2 e^{-(0)}\cos(\pi x) + 2\cos(e^{-(0)}(1-\cos(\pi x)))e^{-(0)}(-1+\cos(\pi x))$$

$$= -a_2 + (a_1)_{xx} - 2\cos(a_0)a_1 +$$

$$\pi^2\cos(\pi x)$$

$$+ 2\cos((1-\cos(\pi x)))(-1+\cos(\pi x))$$

$$= -1 + \cos(\pi x) - \pi^2\cos(\pi x) - 2\cos(1-\cos(\pi x))(-1+\cos(\pi x)) + \pi^2\cos(\pi x) + 2\cos((1-\cos(\pi x)))(-1+\cos(\pi x))$$

$$= -1 + \cos(\pi x),$$

Now, by using the equation (9),

$$w(x,t) = a_0 + a_1 t + a_2 \frac{t^2}{2!} + a_3 \frac{t^3}{3!} + \cdots + a_n \frac{t^n}{n!} + \cdots$$

$$= (1 - \cos(\pi x)) - (1 - \cos(\pi x))t +$$
$$(1 - \cos(\pi x))\frac{t^2}{2!} - (1 - \cos(\pi x))\frac{t^3}{3!} + \cdots \quad (2.12)$$

$$= (1 - \cos(\pi x))(1 - t + \frac{t^2}{2!} - \frac{t^3}{3!} + \cdots)$$

$$= (1 - \cos(\pi x))e^{-t}.$$

Then, the graph of exact, analytical solutions, and absolute errors of example 3 for $t = 1$ are given in **Erro! Fonte de referência não encontrada.**. **Erro! Fonte de referência não encontrada.** shows the absolute error using the proposed technique, ADM, and HAM of example 3 with $t = 1, k = 0.001$. When the value of $n$ increases, the approximate solution gradually approaches the exact solution and the absolute error will decrease.

## 3. RESULTS AND DISCUSSION:

Three test problems are introduced for confirming the validity of the novel proposed technique, which used to find the solutions of linear and nonlinear one-dimensional telegraph equations. Figures (1-3) showed that the exact solution, analytical solution, and absolute errors at $t = 1$, $k = 0.001$. Also, the analytical solutions obtained by a proposed technique have been compared with the solutions obtained by ADM (Biazar and Ebrahimi 2007), and HAM (Raftari and Yildirim 2012) in three test examples, which are given in Tables (1-3). We then found that the analytical solutions obtained by a proposed technique are more identical to the exact solutions than those obtained using ADM and HAM of linear and nonlinear one-dimensional telegraph equation. Therefore, it can be seen that the absolute errors of the proposed technique are smaller than ADM and HAM, which are shown in Tables (1-3). More precisely, the measurement of absolute errors for these examples guarantees the ability of the proposed technique and its accuracy in finding the analytical solutions of linear and nonlinear one-dimensional telegraph equation. Moreover, according to the computations that are introduced in the figures and tables, we can say that, the analytical technique is an effective and good technique to find the solutions of linear and nonlinear one-dimensional telegraph equations compared to the ADM, and HAM.

### 3.1. Convergence analysis

Consider the PDE (1.1) in the following form:

$$w(x,t) = G\big(w(x,t)\big), \qquad (1)$$

where $G$ is a nonlinear operator. The solution that obtained by the presented technique is equivalent to the following sequence:

$$S_n = \sum_{i=0}^{n} w_i = \sum_{i=0}^{n} \delta_i \frac{(\Delta t)^i}{(i)!}. \qquad (2)$$

**Theorem 3.1** (Convergence of linear and nonlinear telegraph equation)

Let $G$ be an operator from a Hilbert space $H$ into $H$ and $w$ be the exact solution of equation (3.1). The approximate solution $\sum_{i=0}^{\infty} w_i = \sum_{i=0}^{\infty} \delta_i \frac{(\Delta t)^i}{(i)!}$ is converged to exact solution $w$, when $\exists\, 0 \leq \delta < 1$, $\|w_{i+1}\| \leq \delta \|w_i\| \forall\, i \in \mathbb{N} \cup \{0\}$.

**Proof:** We want to prove that $\{S_n\}_{n=0}^{\infty}$ is a converged Cauchy sequence,

$$\|S_{n+1} - S_n\| = \|w_{n+1}\| \leq \delta \|w_n\| \leq \delta^2 \|w_{n-1}\| \leq \cdots \leq \delta^n \|w_1\| \leq \delta^{n+1} \|w_0\|. \qquad (3)$$

Now for $n, m \in \mathbb{N}, n \geq m$ we get

$$\|S_n - S_m\| = \|(S_n - S_{n-1}) + (S_{n-1} - S_{n-2}) + \cdots + (S_{m+1} - S_m)\|$$

$$\leq \|S_n - S_{n-1}\| + \|S_{n-1} - S_{n-2}\| + \cdots + \|S_{m+1} - S_m\|$$

$$\leq \delta^n \|w_0\| + \delta^{n-1} \|w_0\| + \cdots + \delta^{m+1} \|w_0\|$$

$$\leq (\delta^{m+1} + \delta^{m+2} + \cdots + \delta^n) \|w_0\| = \delta^{m+1} \frac{1 - \delta^{n-m}}{1 - \delta} \|w_0\|. \qquad (4)$$

Hence, $\lim_{n,m \to \infty} \|S_n - S_m\| = 0$, i. e., $\{S_n\}_{n=0}^{\infty}$ is a Cauchy sequence in the Hilbert space $H$. Thus, there exist $S \in H$ such that $\lim_{n \to \infty} S_n = S$, where $S = w$.

**Definition 3.2** For every $n \in \mathbb{N} \cup \{0\}$, we define

$$\delta_n = \begin{cases} \frac{\|w_{n+1}\|}{\|w_n\|}, & \|w_n\| \neq 0 \\ 0, & otherwise. \end{cases} \qquad (5)$$

**Corollary 3.3** From Theorem 3.1,

$$\sum_{i=0}^{\infty} w_i = \sum_{i=0}^{\infty} \delta_i \frac{(\Delta t)^i}{(i)!},$$

is converged to exact solution $w$ when $0 \leq \delta_i < 1$, $i = 0,1,2,\dots$.

To illustrate the convergence of analytical approximate solutions for the three examples, it was applied Corollary 3.3 as follows.

In the first example where $(x,t) \in (0,1)^2$,

$$\delta_1 = \frac{\|w_2\|}{\|w_1\|} = 0.7745966694 < 1,$$

$$\delta_2 = \frac{\|w_3\|}{\|w_2\|} = 0.5634361702 < 1,$$

$$\delta_3 = \frac{\|w_4\|}{\|w_3\|} = 0.4409585516 < 1.$$

In the second example where $(x,t) \in (0,1)^2$,

$$\delta_1 = \frac{\|w_2\|}{\|w_1\|} = 0.06211299938 < 1,$$

$$\delta_2 = \frac{\|w_3\|}{\|w_2\|} = 0.02773500981 < 1,$$

$$\delta_3 = \frac{\|w_4\|}{\|w_3\|} = 0.01561561843 < 1.$$

In the third example where $x \in (0,2)$ and $t \in (0,1)$,

$$\delta_1 = \frac{\|w_2\|}{\|w_1\|} = 0.3872983346 < 1,$$

$$\delta_2 = \frac{\|w_3\|}{\|w_2\|} = 0.2817180849 < 1,$$

$$\delta_3 = \frac{\|w_4\|}{\|w_3\|} = 0.2204792759 < 1.$$

Therefore, the convergence of analytical solutions is valid. Finally, the theoretical proofs for the analysis of convergence coincide with the computation results presented in the above figures and tables.

## 4. CONCLUSIONS:

The proposed technique is an efficient methodology with good accuracy and convergence and a powerful tool to find approximate analytic solutions for the linear and nonlinear problems. The tests confirm the validity of a novel technique to handle current linear and nonlinear PDEs. In the future, this research can be extended to the investigation by applying this technique for more complicated problems such as systems of nonlinear PDEs.

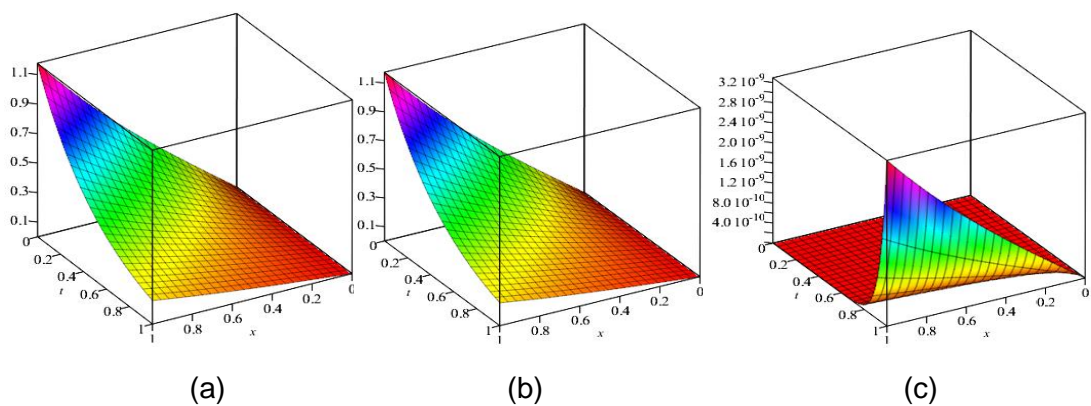

**Figure 1.** *(a) Exact solution, (b) Approximate solution, (c) Absolute errors for example 1 when t= 1.*

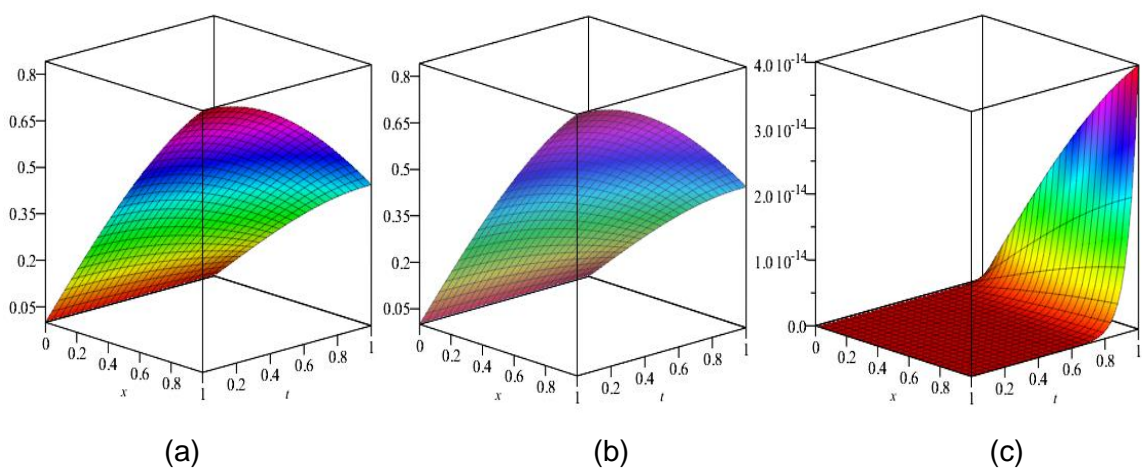

**Figure 2.** *(a) Exact solution, (b) Approximate solution, (c) Absolute errors for example 2 when t= 1.*

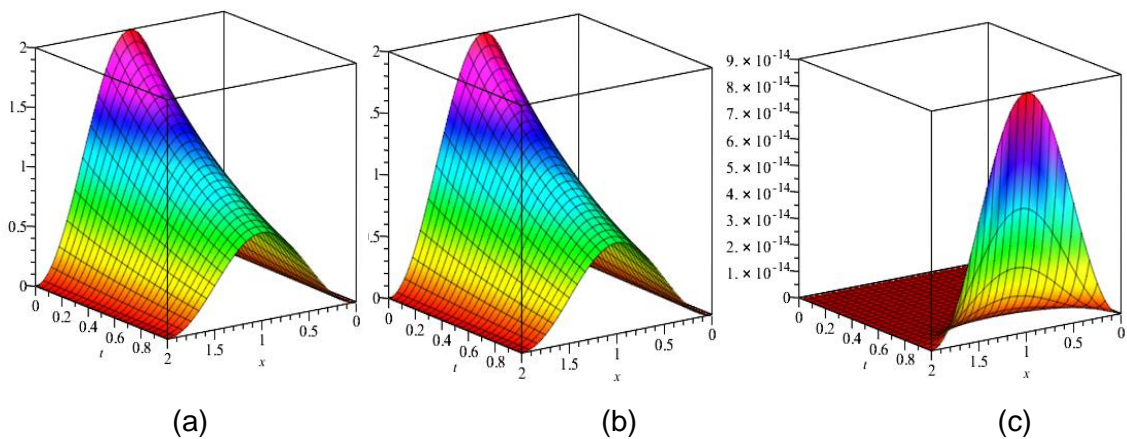

**Figure 3.** *(a) Exact solution, (b) Approximate solution, (c) Absolute errors for example 3 when t= 1.*

**Table 1.** *Comparison of absolute errors among different methods and present method, for example, 1 with $t = 1,\ k = 0.001, \propto = 10, \beta = 5$.*

| Method | x | n=10 | n=11 | n=12 | n=13 | n=14 | n=15 |
|---|---|---|---|---|---|---|---|
| Proposed method | | 8.84E-06 | 1.49E-06 | 2.32E-07 | 3.34E-08 | 4.48E-09 | 5.64E-10 |
| ADE method | 0.1 | 2.44E-05 | 2.87E-04 | 3.24E-06 | 4.42E-06 | 2.71E-06 | 2.16E-08 |
| HAM method | | 1.54E-05 | 3.51E-05 | 1.92E-05 | 2.90E-05 | 3.84E-07 | 6.47E-07 |
| Proposed method | | 1.80E-05 | 3.04E-06 | 4.72E-07 | 6.81E-08 | 9.14E-09 | 1.15E-09 |
| ADE method | 0.3 | 2.23E-05 | 6.14E-05 | 6.11E-05 | 4.76E-06 | 7.66E-08 | 2.16E-06 |
| HAM method | | 3.54E-04 | 4.27E-04 | 1.53E-03 | 5.33E-07 | 3.27E-07 | 3.83E-08 |
| Proposed method | | 2.80E-05 | 4.71E-06 | 7.32E-07 | 1.05E-07 | 1.42E-08 | 1.78E-09 |
| ADE method | 0.5 | 5.11E-04 | 6.03E-04 | 6.03E-06 | 3.88E-06 | 6.56E-05 | 4.11E-05 |
| HAM method | | 3.03E-04 | 3.62E-04 | 4.16E-05 | 7.36E-04 | 2.45E-06 | 2.74E-07 |
| Proposed method | | 3.91E-05 | 6.57E-06 | 1.02E-06 | 1.47E-07 | 1.98E-08 | 2.49E-09 |
| ADE method | 0.7 | 6.54E-04 | 2.82E-05 | 7.15E-03 | 4.11E-05 | 5.78E-07 | 4.33E-06 |
| HAM method | | 4.76E-03 | 5.67E-05 | 5.44E-04 | 7.54E-04 | 2.14E-05 | 2.54E-08 |
| Proposed method | | 5.16E-05 | 8.70E-06 | 1.35E-06 | 1.95E-07 | 2.62E-08 | 3.29E-09 |
| ADE method | 0.9 | 6.23E-04 | 4.65E-03 | 7.73E-04 | 5.22E-05 | 6.25E-07 | 6.15E-07 |
| HAM method | | 3.90E-03 | 3.32E-04 | 4.51E-05 | 3.12E-06 | 5.32E-06 | 3.41E-06 |

**Table 2.** *Comparison of absolute errors among different methods and present method for example 2 with $t = 1,\ k = 0.001, \propto = 10, \beta = 5$.*

| Method | x | n=10 | n=11 | n=12 | n=13 | n=14 | n=15 |
|---|---|---|---|---|---|---|---|
| Proposed method | | 4.12E-10 | 4.12E-10 | 2.27E-12 | 2.27E-12 | 9.46E-15 | 9.46E-15 |
| ADE method | 0.1 | 3.54E-09 | 2.23E-08 | 4.01E-09 | 4.90E-10 | 7.21E-12 | 6.87E-11 |
| HAM method | | 2.04E-07 | 3.64E-07 | 3.66E-10 | 2.11E-11 | 5.99E-11 | 7.33E-09 |
| Proposed method | | 8.09E-10 | 8.09E-10 | 4.45E-12 | 4.45E-12 | 1.86E-14 | 1.86E-14 |
| ADE method | 0.3 | 6.64E-08 | 5.33E-07 | 3.49E-10 | 1.13E-10 | 2.18E-11 | 4.21E-11 |
| HAM method | | 6.72E-08 | 6.97E-09 | 2.19E-10 | 2.91E-11 | 3.02E-10 | 2.99E-12 |
| Proposed method | | 1.17E-09 | 1.17E-09 | 6.45E-12 | 6.45E-12 | 2.69E-14 | 2.69E-14 |
| ADE method | 0.5 | 3.06E-05 | 3.91E-07 | 5.03E-10 | 4.51E-11 | 6.91E-10 | 3.42E-10 |
| HAM method | | 3.77E-07 | 5.64E-06 | 3.16E-10 | 4.18E-10 | 3.02E-11 | 7.11E-11 |
| Proposed method | | 1.49E-09 | 1.49E-09 | 8.19E-12 | 8.19E-12 | 3.42E-14 | 3.42E-14 |
| ADE method | 0.7 | 3.33E-06 | 5.55E-06 | 6.67E-10 | 5.61E-08 | 3.91E-11 | 7.98E-10 |
| HAM method | | 6.54E-06 | 3.96E-06 | 4.88E-09 | 7.98E-10 | 6.22E-10 | 3.90E-12 |
| Proposed method | | 1.75E-09 | 1.75E-09 | 9.61E-12 | 9.61E-12 | 4.01E-14 | 4.01E-14 |
| ADE method | 0.9 | 6.32E-06 | 4.36E-06 | 7.21E-09 | 8.19E-10 | 2.18E-10 | 2.66E-11 |
| HAM method | | 3.55E-06 | 3.28E-06 | 5.27E-10 | 6.03E-09 | 2.80E-09 | 4.18E-10 |

| Method | x | n=10 | n=11 | n=12 | n=13 | n=14 | n=15 |
|--------|---|------|------|------|------|------|------|
| Proposed method | | 1.13E-09 | 9.48E-11 | 7.33E-12 | 5.26E-13 | 3.52E-14 | 2.21E-15 |
| ADE method | 0.1 | 2.04E-07 | 3.47E-10 | 4.02E-10 | 7.04E-11 | 4.41E-10 | 1.49E-12 |
| HAM method | | 2.44E-08 | 5.97E-09 | 3.83E-09 | 6.48E-10 | 3.28E-11 | 4.47E-10 |
| Proposed method | | 9.53E-09 | 7.99E-10 | 6.18E-11 | 4.43E-12 | 2.97E-13 | 1.86E-14 |
| ADE method | 0.3 | 2.42E-07 | 1.72E-08 | 5.24E-08 | 3.04E-10 | 2.72E-11 | 3.16E-10 |
| HAM method | | 3.07E-06 | 7.92E-07 | 2.47E-09 | 2.82E-09 | 6.21E-10 | 2.45E-12 |
| Proposed method | | 2.31E-08 | 1.94E-09 | 1.50E-10 | 1.08E-11 | 7.20E-13 | 4.51E-14 |
| ADE method | 0.5 | 5.11E-05 | 3.64E-07 | 6.12E-07 | 4.63E-10 | 2.57E-09 | 6.77E-08 |
| HAM method | | 4.37E-06 | 5.23E-06 | 4.69E-08 | 6.99E-08 | 6.02E-09 | 4.92E-09 |
| Proposed method | | 3.67E-08 | 3.08E-09 | 2.38E-10 | 1.71E-11 | 1.14E-12 | 7.17E-14 |
| ADE method | 0.7 | 8.20E-06 | 7.87E-06 | 7.74E-06 | 3.42E-09 | 2.11E-09 | 7.52E-11 |
| HAM method | | 4.65E-05 | 4.43E-07 | 4.04E-08 | 5.11E-08 | 7.06E-10 | 3.82E-10 |
| Proposed method | | 4.51E-08 | 3.78E-09 | 2.92E-10 | 2.10E-11 | 1.40E-12 | 8.81E-14 |
| ADE method | 0.9 | 4.06E-07 | 3.87E-07 | 4.21E-08 | 4.01E-06 | 3.93E-09 | 3.01E-10 |
| HAM method | | 6.27E-06 | 5.90E-07 | 5.64E-09 | 1.23E-09 | 6.73E-07 | 4.25E-08 |